\documentclass[11pt,twoside]{article}
\usepackage{a4}
\usepackage{amssymb,amsmath,amsthm,latexsym}
\usepackage{amsfonts}
	\usepackage{amsfonts}
\usepackage{graphicx}
\usepackage[pdftex,bookmarks,colorlinks=false]{hyperref}

\def\ni{\noindent}

\vspace{5cm}

\begin{document}
	
\label{'ubf'}
\setcounter{page}{1} 

\markboth {\hspace*{-9mm} \centerline{\footnotesize \sc
  Point- and arc-reaching sets in a digraph}
                 }
                { \centerline {\footnotesize \sc
         B.D. Acharya, Germina K.A., Kumar Abhishek, S.B. Rao and T. Zaslavsky } \hspace*{-9mm}
               }

\begin{center}

       {\Large \textbf { \sc Point- and arc-reaching sets of vertices in a digraph
                               }
       }
\\[10pt]

{\sc B.D. Acharya$^\S,$}
{\sc K.A. Germina$^\dag,$}
{\sc Kumar Abhishek$^\ddag,$}\\
{\sc S.B. Rao$^\sharp,$} and
{\sc Thomas Zaslavsky$^*.$}
\\[10pt]
{\footnotesize $^\S$ Srinivasa Ramanujan Center for Intensification of Interaction in Interdisciplinary \\Discrete Mathematical Sciences (SRC-IIIDMS),\\ University of Mysore, Mysore - 560 005, India. \footnote{\emph{Address for correspondence}: No.22, $10^{th}$ Cross, $5^{th}$ Main, New Thippasandra Post, Malleshpalya, Bangalore-560 075, India.}\\
{\footnotesize $^{\dag, \ddag}$Research Center \& PG Department of Mathematics, \\ Mary Matha Arts \& Science College, \\ Vemom P.O., Mananthavady - 670645, \ India.}\footnote{$^\ddag$\underline{Present address}: Department of Mathematics, Amrita Vishwa Vidya Peetham, Ettimadai, Coimbatore-641105, India.}\\
{\footnotesize $^\sharp$C.R. Rao Advanced Institute for Mathematics, Statistics and Computer Science, \\ `ARYABHATA', University of Hyderabad Campus, \\ Central University P.O., Hyderabad-500046, India.}\\
{\footnotesize $^*$Department of Mathematical Sciences,\\ Binghamton University (SUNY), \\ Binghamton, NY 13902-6000, U.S.A.}\\[10pt]
{\footnotesize e-mail: {\it devadas.acharya@gmail.com$^\S$}}\\
{\footnotesize e-mail: {\it srgerminaka@gmail.com$^\dag$ }}\\
{\footnotesize e-mail: {\it sethkumarabhishek@hotmail.com$^\ddag$}}\\
{\footnotesize e-mail:{\it siddanib@yahoo.co.in$^\sharp$}}\\
{\footnotesize e-mail:{\it zaslav@math.binghamton.edu$^*$}}

}

\end{center}
\thispagestyle{empty}

\hrulefill

\begin{center} \textbf{Abstract} \end{center}

\small

In a digraph $D = (X, \mathcal{U})$, not necessarily finite, an arc $(x, y) \in \mathcal{U}$ is reachable from a vertex $u$ if there exists a directed walk $W$ that originates from $u$ and contains $(x, y)$. A subset $S \subseteq X$ is an arc-reaching set of $D$ if for every arc $(x, y)$ there exists a diwalk $W$ originating at a vertex $u \in S$ and containing $(x, y)$. A minimal arc-reaching set is an arc-basis. $S$ is a point-reaching set if for every vertex $v$ there exists a diwalk $W$ to $v$ originating at a vertex $u \in S$. A minimal point-reaching set is a point-basis. We extend the results of Harary, Norman, and Cartwright on point-bases in finite digraphs to point- and arc-bases in infinite digraphs.
\\

\ni\underline{Key Words}: Digraph, reachability, point-basis, arc-basis.

\ni\underline{Mathematics Subject Classification (2010)}: Primary 05C20. [Erroneously published as 05C22.]

\hrulefill

\section{Introduction}

D\'enes K\"onig, in \cite{dk}, Section VII.1, introduced the concept of a point-basis in a digraph. Letting a {\it point-reaching set} be a vertex set from which there is a directed path to every vertex, a {\it point-basis} is a minimal point-reaching set. K\"onig characterized point-bases in a weak sense and Harary, Norman, and Cartwright, in \cite{hnc}, Chapter 4, gave a stronger characterization for finite digraphs. We give a complete description of all point-reaching sets, point-bases, and their arc analogs in infinite digraphs.

We refer the reader to Berge \cite{cb} and Harary, Norman, and Cartwright \cite{hnc} for all terminology and notation in the theory of digraphs and graphs that are not defined herein.

\section{Preliminaries}

We begin with definitions of the essential concepts. $D$ is a finite or infinite digraph with {\it vertex} (or {\it point}) set $X(D)$ and {\it arc set} $\mathcal{U}(D)$. An {\it isolate} is a vertex which has out-degree and in-degree 0. $D^0$ is $D$ with all isolates removed. $D^-$ is $D$ with all {\it sinks} (points of out-degree 0) removed. For $u, v \in X(D)$, $v$ is {\it reachable} from $u$ if there is a {\it diwalk} (directed walk), or equivalently a dipath, in $D$ from $u$ to $v$. The {\it reach} of $u$ is the set $$R(u) := \{v \in X(D) : v \ \mbox{is reachable from} \ u\}.$$ \ni For $S \subseteq X(D)$ and $y \in X(D)$, $y$ is reachable from $S$ if it is reachable from some vertex in $S$, equivalently, if it is in the set $$R(S) := \bigcup_{v \in S} R(v).$$

A {\it strong component} of $D$ is a maximal subdigraph in which every vertex is reachable from every other.\\

\ni\textbf{Observation 1.} If $v \in R(u)$, then $R(v) \subseteq R(u)$. If $u, v$ are vertices in the same
strong component of $D$, then $R(u) = R(v)$. Every reach $R(u)$ is a union of strong components of $D$.\\

A {\it point-reaching set} is a subset $S \subseteq X(D)$ such that for every point $u \in X(D)$ there
is a dipath to $u$ from a vertex in $S$. A {\it point-basis} is a minimal point-reaching set. An
{\it arc-reaching set} is a subset $S \subseteq X(D)$ such that for every arc $(u, v)$, there is a dipath to $u$ from a vertex in $S$. An {\it arc-basis} is a minimal arc-reaching set. Note that point- and arc-bases do not necessarily exist in an infinite digraph.\\

If $T \subseteq X(D)$, a $T$-{\it reaching set} is a subset $S \subseteq X(D)$ such that for every point $u \in T$ there is a dipath in $D$, to $u$, from some vertex in $S$; i.e., such that $T \subseteq R(S)$. Thus,
an arc-reaching set is a $T(D)$-reaching set where $T(D)$ is the set of tails of arcs in $D$. A point-reaching set is an $X(D)$-reaching set.\\

\ni\textbf{Observation 2.} The point-reaching sets of $D$, the arc-reaching sets of $D$,
and in general the $T$-reaching sets for any $T \subseteq X(D)$ form {\it superhereditary families of subsets} of $X(D)$; i.e., if $A$ is one such set, and $A \subset B$, then $B$ is also such a set.\\

\ni\textbf{Observation 3.} The arc-reaching sets of $D$ are the arc-reaching sets of $D^0$,
together with any subset of the isolates of $D$. The arc-bases of $D$ are the arc-bases of $D^0$.\\

\ni\textbf{Theorem 4.} An arc-reaching set of a finite or infinite digraph $D$ is identical
to a point-reaching set of $D^-$ together with any subset of the sinks of $D$. An arc-basis of $D$ is identical to a point-basis of $D$ with the
isolates deleted.\\

\emph{Proof.} $X(D^-) = T(D)$. Because a point not in $D^-$ has no out-arcs in $D$, for any two points
$u, v$ in $D^-$, there is a dipath from $u$ to $v$ in $D^-$ if and only if there is such a dipath in $D$. $\Box$\\

Theorem 4 implies that every result about point-reaching sets and point-bases translates into a result about arc-reaching sets and arc-bases. Many such results are stated in \cite{hnc}; here we give only a characterization of point-reaching sets and point-bases in infinite digraphs, Theorem 16, that generalizes the finite case solved in \cite{hnc}, Theorem 4.4.

\section{Characterization of point-bases}

A strong component of a digraph $D$ is {\it initial} if there are no arcs entering it. An isolate, for instance, is an initial strong component (of a particularly simple kind). A {\it source} is a vertex without in-arcs; i.e., it is an initial strong component of order 1 (thus, trivially, every isolate is a source). The {\it condensation} $D^*$ of $D$ is the digraph which has one vertex for each strong component $C$ of $D$, with an arc $(C_1,C_2)$ if and only if there is an arc $(u, v)$ in $D$ such that $u \in C_1$ and $v \in C_2$. The condensation of a vertex set $S \subseteq X(D)$ is the vertex set $$S^* := \{C : C \ \mbox{is a strong component of} \ D \ \mbox{and} \ C \cap S \ne \emptyset\} \subseteq X(D^*).$$ We write $R^*(C)$ for
the set of vertices of $D^*$ (i.e., strong components of $D$) which are reachable from $C$
in $D^*$.\\

\ni\textbf{Observation 5.} A strong component is initial in $D$ if and only if it is a source
vertex in $D^*$.\\

The next observation is \cite{hnc}, Theorem 3.3, but applied to infinite as well as finite graphs.\\

\ni\textbf{Observation 6.} Let $u_i$ be a vertex in the strong component $C_i$ of $D$, for
$i = 1, 2$. Then $u_2$ is reachable from $u_1$ in $D$ if and only if $C_2$ is reachable from
$C_1$ in $D^*$. Furthermore, $R^*(C_1)$ is the condensation of $R(u_1)$.\\

The {\it shadow} of $u \in X(D)$ is $R^{\prime}(u) := \{z \in X(D) : u \ \mbox{is reachable from} \ z\}$.\\

\ni\textbf{Observation 7.} The shadow $R^{\prime}(u)$ equals the reach of $u$ in the converse digraph
$D^{\prime}$ of $D$. Every shadow $R^{\prime}(u)$ is a union of strong components of $D^{\prime}$.\\

An {\it in-ray} in $D$ is a one-way infinite dipath $(. . . , v_k, . . . , v_1, v_0)$ whose vertex set is not contained in the union of finitely many strong components of $D$.\\

\ni\textbf{Example 8.} Consider two digraphs with vertex set $Y := \{y_0, y_1, y_2, . . .\}$ defined as follows. The arcs of $D$ are $(y_{i+1}, y_i)$ for $i = 0, 1, . . .$. In the converse $D^{\prime}$ the arcs are $(y_i, y_{i+1})$ for $i = 0, 1, . . .$. In $D$ and its converse $D^{\prime}$, a set is an arc-basis if and only if it is a point-basis. In $D^{\prime}$ there is one point-basis, $\{y_0\}$. In $D$ there is no point-basis. Any point-reaching set must contain vertices $y_i$ with arbitrarily large subscripts. Thus, every point-reaching set is infinite, and furthermore there is an infinite descending chain of point-reaching sets, for instance the sets $Y_k := \{y_i : i \ge k\}$, whose intersection is empty. (That is about as far as one can get from the existence of a minimal point-reaching set.)\\

\ni\textbf{Example 9.} Let $D$ be the infinite binary tree, in which $X(D)$ is the set of all finite binary
sequences (sequences of 0s and 1s, including the null sequence) and there is an arc $(b_1b_2 \dots b_lb_{l+1}, b_1b_2 \dots b_l)$ for every binary sequence $b_1b_2 \dots b_lb_{l+1}$ with $l \ge 0$. This is a 'more infinite' version of $D$ in Example 8. As in that example, all point-reaching sets are arc-reaching sets, all such sets are infinite, and there is no point-basis or arc-basis.\\

\ni\textbf{Example 10.} Let $D$ be the digraph with vertex set $Y \cup Z$ where $Y := \{y_0, y_1, y_2, . . .\}$ and $Z := \{z_0, z_1, z_2, . . .\}$ and with the arc set $\{(y_i, z_i), (z_i, y_i), (y_{i+1}, y_i), (z_{i+1}, z_i) : i = 0, 1, . . .\}$. The strong components of $D$ are the sets $C_i := \{y_i, z_i\}$. Two of the in-rays are $(. . . , y_2, y_1, y_0)$ and $(. . . , z_2, z_1, z_0)$. The condensation $D^*$ is an in-ray, $(. . . ,C_2,C_1,C_0)$; and $D^- = D$. The point-reaching sets in $D$ are the same as the arc-reaching sets, by Theorem 4. An example of such a set is $Y$. Indeed, a subset of $Y \cup Z$ is a point-reaching set if and only if it is infinite. Hence, $D$ has no point-bases or arc-bases.
The complement of $Y$ contains the in-ray $(z_0, z_1, z_2, . . .)$. Thus, it is possible for a digraph
to have a point-reaching set and an in-ray that contains no member of the point-reaching set.\\

\ni\textbf{Example 11.} Let $D$ be the digraph with vertex set $Y \cup Z$, as in Example 10,
and arc set $\{(y_{i+1}, y_i), (z_{i+1}, y_i) : i = 0, 1, . . .\}$. The condensation $D^* = D$, and $D^- = D$. The in-rays of this digraph are $(. . . , y_{k+2}, y_{k+1}, y_k)$ for each $k = 0, 1, 2, . . .$. The unique point-basis is $Z$, which is also the unique arc-basis. The point-reaching and arc-reaching sets
are precisely the sets that contain $Z$. The complement of the unique point-basis is an in-ray,
$(. . . , y_2, y_1, y_0)$. In this digraph, not only is there an in-ray which contains no member of
the point-basis, but indeed, no strong component that meets the in-ray contains a member of the point-basis.\\

\ni\textbf{Example 12.} Let $D$ be the digraph whose vertex set is $\{u\}\cup U_1 \cup U_2 \cup \cdots$ where $U_i := \{u_{i0}, u_{i1}, . . . , u_{ii}\}$, and whose arc set is $\{(u_{i0}, u) : i = 1, 2, . . .\} \cup \{(u_{ik}, u_{i(k-1)}) : i = 1, 2, . . . ,\ 1 \le k \le i\}$. This infinite digraph contains no in-rays. Its sole point-basis (also its sole arc-basis) is the set $\{u_{ii} : i = 1, 2, . . .\}$.\\

\ni\textbf{Lemma 13.} For any strong component $C$ of $D$, there is an initial strong
component $C_0$ such that $C$ is reachable from $C_0$, or there is an in-ray that ends in $C$, or
both. In the former case, $C_0$ is a source in $D^*$ and $C \in R^*(C_0)$. In the latter case, there is
an in-ray in $D^*$ that ends at $C$.\\

\emph{Proof.} Assume that $C$ is not reachable from any initial strong component. Equivalently (by
Observation 5), $R^{*\prime}(C)$ contains no source. Because $C$ is not a source, there is a vertex
$C_1 \in R^{*\prime}(C)$ such that $(C_1,C)$ is an arc in $D^*$. Since $C_1$ is not a source, there is a vertex $C_2 \in R^{*\prime}(C_1)$ such that $(C_2,C_1)$ is an arc in $D^*$. Continuing to track backward in this fashion we find a sequence of vertices of $D^*$, \ $\{C,C_1,C_2, . . .\}$, such that $W^* := (. . . ,C_2,C_1,C)$ is an in-ray to $C$ in $D^*$.

Now we must convert $W^*$ to an in-ray in $D$.  Let $(u_1,v)$ and $(u_{i+1},v_i)$ for $i=1,2, . . .$ be arcs in $D$ such that $v \in C$, $u_i \in C_i$, and $v_i \in C_i$.  Let $P_i$ be a dipath in $C_i$ from $v_i$ to $u_i$.  Then $W := \cdots P_2 P_1$ is an in-ray in $D$, ending at $v \in C$, whose condensation is $W^*$.
$\Box$\\

\ni\textbf{Theorem 14.} Let $D$ be a finite or infinite digraph and let $A \subseteq X(D)$.\\
\textbf{(a)} The following properties of $A$ are equivalent.\\
\hspace*{1.4em} (i) $A$ is a point-reaching set of $D$.\\
\hspace*{1.4em} (ii) $A^*$ is a point-reaching set of $D^*$.\\
\hspace*{1.4em} (iii) $A$ consists of at least one vertex $u$ from each initial strong component of $D$ together with at least one vertex from every shadow $R^{\prime}(u)$ that contains no initial strong components.\\
\hspace*{1.4em} (iv) $A$ consists of at least one vertex $u$ from each initial strong component of $D$ together with an infinite number of vertices in every shadow $R^{\prime}(u)$ that contains no initial
strong components.\\
\textbf{(b)} A point-basis, if one exists, consists of exactly one vertex from each initial strong component of $D$.\\

\emph{Proof.} First, we note that every point-reaching set $A$ contains at least one vertex from each
initial strong component $C$. Suppose $A$ did not contain a vertex from $C$. Since each $u \in A$ is
not in $C$, and there is no dipath from $u$ to $C$ because $C$ is initial, $A$ cannot be a point-reaching
set.

Next, we prove the equivalence of the four statements in \textbf{(a)}. Equivalence of (i) and (ii) is
a consequence of the last part of Observation 6. It is obvious that (iv) implies (iii) and (iii)
implies (i).

To prove that (i) implies (iv), let $A$ be a point-reaching set and $u$ be a vertex whose shadow
contains no initial strong component. Because $A$ is point-reaching, $A \cap R^{\prime}(u) \ne \emptyset$. Recall that $A^*$ is the condensation of $A$ in $D^*$. Let $C_0$ be the strong component that contains $u$, so that $R^{*\prime}(C_0)$ is the shadow of $C_0$ in $D^*$. This shadow is infinite, because it does not contain a source (by Observation 5). For each strong component $C \in R^*(C_0)$ let $P(C)$ be the unique dipath from $C$ to $C_0$ in $D^*$. The set $R^*(A^*)$ of all vertices in $R^*(C_0)$ that are reachable from $A^*$ is equal to $\bigcup_{i=1}^k P(C_i)$. Because $A^*$ is point-reaching in $D^*$, $R^*(C_0) \subseteq R^*(A^*)$; hence, $R^*(A^*)$ is infinite. However, each path $P(C_i)$ is finite, so if $A^* \cap R^{*\prime}(C_0)$ is finite, then $$\bigcup_{C\in A^* \cap R^{*\prime}(C_0)} P(C)$$ \ni is finite. This is a contradiction; therefore, $A^* \cap R^{*\prime}(C_0)$ is infinite so $A \cap R^{\prime}(u)$ is infinite.

Now we characterize point-bases. A point-basis $B$ cannot contain more than one vertex from $C$, because by Observation 1 having a second vertex from $C$ included in $B$ does not increase the set of reachable vertices, so $B$ would not be minimally point-reaching. Furthermore, $B$ cannot contain a vertex from a non-initial strong component. Suppose it did contain a vertex $u$ from a non-initial strong component $C$. Since $C$ is not initial, there is a strong component $C_0$ such that $(C_0,C)$ is an arc in $D^*$. There must be vertices $v \in C_0$ and $z \in B$ such that $v \in R(z)$. Since $u \in R(v)$, it follows that $u \in R(z)$ by Observation 1. Therefore, $R(B - u) = R(B)$, so $B - u$ is a point-reaching set,  contradicting the assumption that $B$ was minimal. This completes the description of point-bases.   $\Box$\\

Note that the second part of Theorem 16 is the generalization of \cite{hnc}, Theorem 4.4, for infinite digraphs.\\

\ni\textbf{Theorem 15.} Let $D$ be a finite or infinite digraph and let $A \subseteq X(D)$. \\
\textbf{(a)} The following properties of $A$ are equivalent.\\
\hspace*{1.4em} (i) $A$ is an arc-reaching set of $D$.\\
\hspace*{1.4em} (ii) $A^*$ is an arc-reaching set of $D^*$.\\
\hspace*{1.4em} (iii) $A$ contains at least one vertex from each initial strong component of $D$ that is not an isolate, as well as at least one vertex from every shadow $R^{\prime}(u)$ that contains no initial strong components.\\
\hspace*{1.4em} (iv) $A$ contains at least one vertex from each initial strong component of $D$ that is not an isolate, as well as an infinite number of vertices in every shadow $R^{\prime}(u)$ that contains no initial strong components.\\
\textbf{(b)} An arc-basis, if one exists, consists of exactly one vertex, not a sink, from each initial strong component of $D$ that is not an isolate.\\

\emph{Proof.} By Theorem 4 this is a corollary of Theorem 14. $\Box$\\

\ni\textbf{Theorem 16.} In a finite or infinite digraph $D$,\\
\textbf{(a)} the following statements are equivalent:\\
\hspace*{1.4em} (i) $D$ has a point-basis.\\
\hspace*{1.4em} (ii) $D$ has an arc-basis.\\
\hspace*{1.4em} (iii) The shadow of every vertex in $D$ contains an initial strong component.\\
\hspace*{1.4em} (iv) The shadow of every vertex in $D^*$ contains a source.\\
\hspace*{1.4em} (v) Every point-reaching set contains a point-basis.\\
\hspace*{1.4em} (vi) Every arc-reaching set contains an arc-basis.\\
\textbf{(b)} A point-basis and an arc-basis always exist if $D$ is finite, more generally if $D$ has a finite number of strong components, and even more generally if $D^*$ does not contain an in-ray.\\

\emph{Proof.} Excepting conditions (v) and (vi), this is an immediate consequence of the description of point- and arc-bases in Theorems 14 and 15 and of Lemma 13.

Lastly, we prove that (i) implies (v) and (ii) implies (vi). (The converses are obvious.)
Assume (i) is true, so that (iii) is also true by the preceding part of the proof. Suppose $A$ is
a point-reaching set. By Theorem 14, $A$ contains a point from each initial strong component
of $D$. Define $A_0$ to consist of one vertex of $A$ from each such component. Let $v \in A$ be such
that $v$ is not in an initial strong component of $D$. Then by (iii), there is an initial strong
component in the shadow of $v$, and by definition this component contains a vertex of $A_0$.
Therefore, $v \in R(A_0)$. It follows that $A \subseteq R(A_0)$, and since $A$ is point-reaching, $A_0$ also is point-reaching. By part (b), $A_0$ is a point-basis. \\

We now deduce that (ii) implies (vi). Since $D$ has an arc-basis, which is a point-basis of
$D^-$ by Theorem 4, $D^-$ has a point-basis. Thus, (i) holds for $D^-$. Consequently, (v) holds
for $D^-$. An arc-reaching set $A$ of $D$ contains a point-reaching set $A^-$ of $D^-$ by Theorem 4,
which contains a point-basis of $D^-$ by (v), which is an arc-basis of $D$ by Theorem 4. $\Box$\\

\ni\textbf{Corollary 17.} Every point-basis in $D$ has the same cardinality. Every
arc-basis in $D$ has the same cardinality.\\

\emph{Proof.} The statement is trivial if there is no point- or arc-basis. Otherwise, the number of
vertices in a point-basis equals the number of initial strong components. The number in an
arc-basis equals the number of initial strong components that are not isolates. $\Box$\\

\ni\textbf{Corollary 18.} Every source lies in every point-generating set. Every source that
is not an isolate lies in every arc-reaching set.

The set of sources is a point-reaching set if and only if every initial strong component has
order 1 and $D$ has no in-rays. It is an arc-reaching set if and only if every initial strong
component has order 1 and $D$ has no isolates and no in-rays.\\

\emph{Proof.} A source is a one-vertex strong component and is initial. The set of sources contains a vertex from every initial strong component if and only if every such component has order 1. $\Box$\\

\ni\textbf{Corollary 19.} Let $A$ be a point-basis of $D$. Then $X(D) - A$ contains a
point-reaching set of $D$ if and only if $D$ has no sources.

Let $A$ be an arc-basis of $D$. Then $X(D)-A$ contains an arc-reaching set of $D$ if and only
if $D$ has no sources other than isolates.

\emph{Proof.} If $D$ has a source $u$, then $u \in A$ so $u \in X(D)-A$. Thus, $u$ is not reachable from any subset of $X(D) - A$.

Because $D$ has a point-basis, every shadow contains an initial strong component, by Theorem 16.
If $D$ has no sources, then every initial strong component contains a vertex that does not belong to $A$.  The set of all such vertices contains a vertex in every initial strong component; thus, it is a point-reaching set by Theorem 14. Hence, $X(D)-A$ contains a point-reaching set of $D$. The same argument applied to $D^-$ proves that $X(D)-A$ contains an arc-reaching set of $D$. $\Box$\\

\ni\textbf{Corollary 20.} Every singleton subset in a digraph is an arc-basis if and only if
the digraph is strongly connected and is not a single point.\\

\emph{Proof.} This is immediate from Theorem 16. $\Box$

\section{A link with domination?}

Curiously, the following result seems to indicate a possibility to generalize a classical theorem on domination in locally finite graphs due to O. Ore~\cite{oo}.\\

\ni\textbf{Theorem 21.} Let $D = (X,\ \mathcal{U})$ be any digraph having at least one arc and such that every infinite descending chain of arc-reaching sets in $D$ has nonempty intersection. Then, every arc-reaching set in $D$ contains a minimal one.\\

\emph{Proof.} Suppose that the theorem is false. Then, there exists a digraph $D$ having an arc-reaching set $S_1$ that contains no minimal one. Since $S_{1}$ does not contain any arc-basis, it in particular is not an arc-basis of $D$, whence there exists an arc-reaching set $S_{2} \subset S_{1}$. Next, by the same argument, since $S_{2}$ is not an arc-basis of $D$ there exists an arc-reaching set $S_{3}\subset S_{2}$. Continuing in this manner, we get an infinite descending chain $S_1 \supset S_2 \supset S_3 \supset \cdots$ of arc-reaching sets of $D$ whose intersection is empty. This is contrary to the hypothesis. Thus, the theorem is established by contraposition.  $\Box$\\

Note that Theorem 16 (i, v) and (ii, vi) imply that if a digraph $D$
possesses a point- (arc-)basis, then every point- (arc-)reaching set
of $D$ contains a minimal one.  However, Examples 9 and 10 indicate that
an infinite digraph may have point- (arc-)reaching sets, yet it may
not possess any point- (arc-)bases.  Theorem 16 (iii, iv) characterize
digraphs that do have a point- (arc-)basis, but the characterizations
are not as simple as we would like.  Hence, we believe the following
problem remains open.

\ni\textbf{Problem 22.} Characterize in a simple way an infinite digraph that
possesses a point- (arc-)basis.

\section{Postscript}

We have left another open problem:

\ni\textbf{Problem 23.}  Generalize our results to $T$-reaching sets and minimal
$T$-reaching sets ('$T$-bases'), where $T$ is an arbitrary subset of the vertex set of $D$, or of its
arc set, or of their union.

\vspace{0.80cm}
\centerline{\textbf{Acknowledgements}}

The second and the third authors are thankful to the Department of Science and Technology (DST), Government of India for supporting their research under project No.SR/S4/MS:277/06. The fourth author is thankful to DST for their support to research under DST-CMS project No.SR/S4/MS:516/07.

\end{document}